\documentclass[11pt]{amsart}
\usepackage{verbatim, amscd, epsfig,amsmath,amsthm,amstext,amssymb,hyperref
}

\newcommand{\stern}{\marginpar{$\bigstar$}}
\newcommand{\kommentar}[1]{\stern\marginpar{\tiny  #1}}

\theoremstyle{plain}
\newtheorem{theorem}{Theorem}[section]
\newtheorem{lemma}[theorem]{Lemma}
\newtheorem{cor}[theorem]{Corollary}

\theoremstyle{definition}
\newtheorem{definition}[theorem]{Definition}

\theoremstyle{remark}
\newtheorem{rem}[theorem]{Remark}

\numberwithin{equation}{section}

\renewcommand{\Im}{{\rm Im}\,}
\newcommand{\del}{\partial}
\newcommand{\delbar}{\bar{\del}}

\newcommand{\R}{\mathbb{ R}}
\newcommand{\C}{\mathbb{ C}}

\newcommand{\Zz}{\mathcal{Z}}

\renewcommand{\H}{\mathbb{ H}}

\renewcommand{\P}{\mathbb{ P}}
\newcommand{\HP}{\H\P}
\newcommand{\CP}{\C\P}

\newcommand{\trivial}[1]{\underline{\H}^{#1}}
\newcommand{\dbar}{{\bar{\partial}}}
\newcommand{\invers}{^{-1}}
\DeclareMathOperator{\End}{End}
\DeclareMathOperator{\Hom}{Hom}
\DeclareMathOperator{\pr}{pr}
\DeclareMathOperator{\ord}{ord}
\DeclareMathOperator{\Gl}{GL}
\DeclareMathOperator{\tr}{tr}

\newcommand{\Rr}{{\mathcal R}}

\begin{document}

\title {Envelopes and osculates  of Willmore surfaces}
\author{K. Leschke and F. Pedit}
\address{
Department of Mathematics \\
University of Massachusetts\\
Amherst, MA 01003, USA\\
}
\email{leschke@gang.umass.edu, franz@gang.umass.edu}
\thanks{MSC--class: 53Axx, 53Cxx, 30Fxx\\
Authors partially supported by SFB288 and first author 
  partially supported by  NSF-grant DMS-9626804}
\maketitle

\parindent0em

\section{Introduction}

Special surfaces in $3$ and $4$--space allow transformations
preserving their special properties: classically these are the
B\"acklund transformations of surfaces of constant curvature and the
Darboux transformations of isothermic surfaces and, more recently,
also the B\"acklund and Darboux transformations of Willmore surfaces
\cite{coimbra}.  These transformations allow to construct more
complicated examples of these special surfaces from simple known
examples.

On the other hand, a complex holomorphic curve in $\C\P^n$ gives rise
to new holomorphic curves: the higher osculating curves and envelopes.
Since every conformal surface in $3$ or $4$--space is a (quaternionic)
holomorphic map $f: M\to \H\P^1 = S^4$ of a Riemann surface $M$, we
expect that the geometric constructions of envelopes and osculating
curves for (quaternionic) holomorphic curves in $\H\P^n$ should relate
to the classical transformation theory of special surfaces. It is in
this spirit that we study osculating and enveloping constructions for
holomorphic curves $f:M\to \H\P^n$ and their effects on conformally
parameterized surfaces with special emphasis on Willmore surfaces.

One of the obstacles in defining osculating curves for a holomorphic
curve $f:M\to \H\P^n$ lies in the fact that the osculating flag
\[
L\subset V_1\subset\dots \subset V_{n-1}\subset V=M \times \H^{n+1}
\] 
built of the successive derivatives of $f$ is only continuous
\cite{Klassiker} along the Weierstrass points $D\subset M$ of $f$.
This motivates the study of a more restricted class of holomorphic
curves $f:M\to\H\P^n$, the {\em Frenet curves}. For these curves the
osculating flag is smooth on $M$ and so is the so--called
\emph{canonical complex structure} $S$ on the trivial
$\trivial{n+1}$--bundle $V$ over $M$. This complex structure is the
analog of the mean curvature sphere congruence along a conformal
immersion into 3 or 4--space.  In particular, $S$ stabilizes the
osculating flag $V_k$ and renders all $V_k$ holomorphic as maps into
the appropriate quaternionic Grassmannians.

We show that the class of Frenet curves is closed under osculating and
enveloping constructions: the \emph{$k$-th osculating curve} $f_k$ of
a Frenet curve $f$ is obtained by intersecting $V_k$ with a
$\H\P^{n-k}\subset \H\P^n$ and is again Frenet. The \emph{envelope} of
a Frenet curve $\tilde{f}:M\to \H\P^{n-1}$ is a holomorphic curve
$f:M\to \H\P^n$ whose first osculating curve $f_1$ equals $\tilde{f}$.
Provided that the space of tangents $H^0(K\tilde{L})$ has a nowhere
vanishing section $\omega$ there is an envelope $f$ of $\tilde{f}$
satisfying $df=\omega$ which is again Frenet.  Since the components
$f=[f_0:\dots :f_n]:M\to \H\P^n$ of a holomorphic curve are conformal
maps into 4-space, the osculating and enveloping constructions build
families of new conformal maps of a Riemann surface $M$ from a given
one by differentiation, integration and algebraic manipulations.

Finally, we define the notion of a Willmore curve into $\HP^n$ and
discuss their osculating and enveloping constructions: Willmore curves
are Frenet curves $f:M\to\H\P^n$ which are critical for the Willmore
energy with respect to variations $f_t:M_t\to \H\P^n$ by Frenet
curves, including variations of the complex structure on $M$. For
$n=1$ these are the usual Willmore surfaces in 3 and 4--space
\cite{willmore_book}.  The Euler-Lagrange equation expresses the fact
that the canonical complex structure $S$ on $V$ is harmonic.
It turns out that, at least for Willmore spheres, the geometric
constructions of envelopes and osculates preserve the Willmore
property. This allows us to construct families of Willmore spheres in
3 and 4--space from a given Willmore sphere by differentiation,
integration and algebraic manipulations. Even if the Willmore sphere
one starts with is a twistor projection of a rational curve into
$\C\P^{2n+1}$, the osculating and enveloping spheres will generally 
not come from rational curves. In particular, we obtain minimal
spheres with planar ends in $\R^4$ via these constructions.

For higher genus Willmore surfaces the osculating and enveloping
constructions need not preserve the Willmore property. Rather,
they give examples of the larger class of {\em constrained} Willmore
surfaces, i.e., holomorphic curves critical for the Willmore energy
under variations fixing the Riemann surface $M$. Of course, in the
case $M=S^2$ constrained Willmore is the same as Willmore. 

\section{Osculates and envelopes  of Frenet curves in $\HP^n$}
\label{sec:frenet}

The successive higher derivatives of a holomorphic curve in $\CP^n$
form a holomorphic flag, the Frenet flag. The intersection of the
$k^{\text{th}}$ osculating flag with a complementary $\CP^{n-k}$ gives
\cite[Ch.\ 2.4]{griffith_harris}
a new holomorphic curve in $\CP^{n-k}$. The analogous construction for
a holomorphic curve in $\HP^n$ requires the existence of a smooth
osculating flag.  From previous work, \cite{Klassiker}, it
is known that a holomorphic curve in $\HP^n$ has a smooth Frenet flag
away from its Weierstrass points into which the flag generally extends
only continuously.

We will briefly recall notions and results of quaternionic holomorphic
geometry, for more details see \cite[Sec.\ 2.5, 4.1, and
4.2]{Klassiker} and \cite[Ch.\ 5 and 6]{coimbra}.

Recall that a map $f$ into $\HP^n$ is the same as a quaternionic line
subbundle $L$ of the trivial $\H^{n+1}$--bundle $V=\trivial{n+1}$ over
$M$, namely $L_p = f(p)$ for $p\in M$.  A smooth map $f: M \to \HP^n$
of a Riemann surface $M$ is \emph{holomorphic}, \cite{Klassiker}, if
there exists a complex structure $J\in \Gamma(\End L)$, $ J^2 =-1$,
such that
\begin{equation}
\label{eq:holo_curve} 
*\delta_L = \delta_L J.
\end{equation}
Here $\delta_L = \pi_{L}\nabla|_L\in \Omega^1(\Hom(L,V/L))$ is the
derivative of $f$ where $\pi_L: V \to V/L$ is the canonical
projection, $\nabla$ is the trivial connection on $V$,
and $*\omega(X) = \omega(J_MX)$ for a $1$-form $\omega$ on $M$.

We say that a holomorphic curve $f: M \to \HP^n$ admits a \emph{Frenet
  flag} if there exists a smooth flag $ V_0 = L \subset V_1 \subset
\ldots \subset V_{n-1} \subset V_n = V $ of quaternionic subbundles of
rank $V_k = k+1$ such that
\begin{equation}
\label{eq:frenet_curve0}
\nabla\Gamma(V_k)\subset \Omega^1(V_{k+1}),
\end{equation}
and a smooth complex structure $S\in\Gamma(\End(V)), S^2=-1$,
stabilizing the flag, with
\begin{equation}
\label{eq:frenet_curve}
 *\delta_{k} = S\delta_{k} = \delta_{k} S.
\end{equation}
As above, $\delta_{k} =
\pi_{V_k}\nabla|_{V_k}\in\Omega^1(\Hom(V_k/V_{k-1},V_{k+1}/V_k))$ are
the derivatives of the $V_k$, where $\pi_{V_k}: V \to V/V_k$ are the
canonical projections.  Note that such an osculating flag is
necessarily unique, whereas the complex structure $S$ is not.
Moreover, $*\delta_{0} = \delta_{0} S$ implies by
\eqref{eq:holo_curve} that $S|_L = J$.

If $f: M \to \HP^n$ is a holomorphic curve then it is shown in
\cite[Lemma 4.1]{Klassiker}  that $f$ admits a Frenet flag away
from a discrete set of points, the \emph{Weierstrass points} $D\subset
M$. These consist of the zeros of the flag derivatives $\delta_k$ as
in the case of complex curves \cite[Ch.\ 2.4]{griffith_harris}. 
Moreover, over $M \setminus D$ there is a unique complex structure of
$V$, the \emph{canonical complex structure}, satisfying
(\ref{eq:frenet_curve}) and additionally the following condition: if
\begin{equation}
\label{eq:nabla_S}
  \nabla S = 2(*Q -*A)
\end{equation}
is the type decomposition of the derivative of a complex structure
$S$, i.e., $Q \in \Gamma(\bar K \End(V))$ and $A\in \Gamma(K\End(V))$,
then
\begin{equation}
\label{eq:canonical_S}
Q|_{V_{n-1}} =0 \quad \text{or, equivalently,} \quad
AV \subset L\,.
\end{equation}
Whereas the Frenet flag generally extends continuously \cite[Lemma
4.10]{Klassiker} across the Weierstrass points $D$, the canonical
complex structure $S$ may become singular as the following example
shows \cite{paule}: a holomorphic curve $f: M \to\HP^1$ is a branched
conformal immersion into $S^4=\HP^1$ whose branch points are the
Weierstrass points $D\subset M$. The Frenet flag $L\subset V$ is
clearly smooth on $M$.  The canonical complex structure $S$ on $V$ is
the mean curvature sphere congruence along $f$ at the immersed points
$p\in M$: the set of eigenlines for $S_p$ on $\H^2$ is a round
2--sphere in $S^4$ which touches $f$ at $p$ by
(\ref{eq:frenet_curve}), and condition (\ref{eq:canonical_S}) says
\cite[Thm.\ 2]{coimbra} that this 2--sphere is the mean curvature sphere
of $f$ at $p\in M$.  If $f: M \to\HP^1$ is the twistor projection of a
complex holomorphic curve $h: M \to\CP^3$ then the mean curvature
sphere is given by the tangent line $W_1\subset V$ of $h$, namely
$S|_{W_1} = i$ and $S|_{W_1j} = -i$. But the tangent $W_1\subset V$ of
$h$ can become quaternionic, i.e., $W_1 = W_1 j$ at some $p\in M$. In
this case the mean curvature sphere $S$ degenerates to a point at
$p\in M$ and thus the complex structure $S$ cannot be extended into
$p\in M$. To avoid these difficulties, we will only consider
holomorphic curves $f: M \to\HP^n$ which have a smooth canonical
complex structure. For conformal maps $f: M \to \HP^1$ this means that
the mean curvature sphere congruence extends smoothly across the
branch points.

\begin{lemma}
\label{lem:S_smooth}
 Let $f: M \to\HP^n$ be  a holomorphic curve with smooth
  canonical complex structure $S$. Then the Frenet flag of $f$ extends
  smoothly across the Weierstrass points.
\end{lemma}
\begin{proof}
  Since
\[ 
*\delta_0 = S\delta_0 = \delta_0 S
\]
the derivative $\delta_0\in H^0(K\Hom_+(L, V/L))$ of $f$ is a complex
holomorphic bundle map \cite{Le}. Here and throughout the paper
$\Hom_\pm$ denote the complex linear respectively complex antilinear
homomorphisms.  Thus $\Im \delta$ defines a smooth quaternionic line
subbundle in $V/L$ whose lift under the canonical projection $\pi:V
\to V/L$ gives the first osculating bundle $V_1\subset V$.  Proceeding
inductively, we extend all the osculating bundles $V_k$ smoothly
across the Weierstrass points.
\end{proof}

\begin{definition} 
A holomorphic curve $f: M \to \HP^n$ is called a \emph{Frenet curve}
if the canonical complex structure, and hence also the Frenet flag,
extends smoothly across the Weierstrass points.
\end{definition}

To construct a first osculating or tangent curve of a Frenet curve $f:
M \to \HP^n$, we choose a hyperplane $H\subset \H^{n+1}$ and intersect
the first Frenet flag $V_1$ with $H$.  If $f$ does not intersect the
hyperplane $H$, we obtain a smooth curve $\tilde f: M \to \HP^{n-1} =
\P(H)$.  By transversality there are hyperplanes not intersecting the
curve.

\begin{definition} The intersection of the first flag $V_1$ of a
  Frenet curve $f: M \to\HP^n$ with a hyperplane $H\subset \H^{n+1}$
  is called the \emph{first osculating} or \emph{tangent curve} of $f$
  with respect to $H$. Conversely, an \emph{envelope} of $f$ is a
  holomorphic curve $\tilde f: M \to\HP^{n+1}$
  whose tangent curve  equals $f$.
\end{definition}
We now show that these constructions preserve Frenet curves.
\begin{theorem}
  A tangent curve of a Frenet curve is Frenet.
\end{theorem}
\begin{proof}
  Let $f: M \to \HP^n$ be a Frenet curve with canonical complex
  structure $S$ and $H\subset \H^{n+1}$ a hyperplane not intersecting
  $f$.  Then $V= L \oplus H$ and $\tilde L = H \cap V_1$ is the
  tangent curve of $f$ so that $V_1 = L \oplus \tilde L$.  Via these
  splittings we identify $H = V/L$ and $\tilde L = V_1/L$.  The flag
  $\tilde V_k = V_{k+1}\cap H$ is the Frenet flag of $\tilde L$ since
  $ \tilde \delta_k = \delta_{k+1} \text{ and } \tilde S = S_{V/L}.$
  Moreover, $\tilde S$ is the canonical complex structure of $\tilde
  L$ since
\[
2*\tilde Q|_{\tilde V_{n-2}} =(\tilde \nabla \tilde
  S)''|_{\tilde V_{n-2}} = \pi(\nabla S)''|_{V_{n-1}\cap H} =
 2 \pi *Q|_{V_{n-1}\cap H} =0\,,
\]
where $\tilde\nabla$ is the connection on $H = V/L$ given by
$\tilde\nabla = \pi \nabla|_H$.
\end{proof}

Applying the theorem successively we obtain 
\begin{cor}
\label{c:Grothendieck}
Let $f: M \to \HP^n$ be a Frenet curve. Then $V$ splits into a direct
sum $V = \oplus_{i=0}^n L_i$ of Frenet curves $f_i : M \to\HP^{n-i}$
where $f_0 = f$ and $f_n$ is a point in $\HP^1$.
\end{cor}

\begin{rem}
  For a Frenet curve $f: M \to \HP^n$ the $k^{\text{th}}$ osculating
  curve $f_k : M \to \HP^{n-k}$, where $L_k = H_k \cap V_k$, depends
  only on the choice of the complementary $n+1-k$--plane $H_k$.
  Therefore one obtains a $k(n+1-k)$--dimensional family of Frenet
  curves in $\HP^{n-k}$. Since the $n-1^{\text{st}}$ osculating curve
  is a map into $\HP^1$, we get a $2(n-1)$--dimensional family of
  branched conformal immersions $f_{n-1}: M \to \HP^1$ via the tangent
  construction.
\end{rem}

To construct an envelope $f$ to a given tangent curve $\tilde f$, we
have to prescribe the possible tangents which $f$ should envelop.
These are the holomorphic $\tilde L$--valued 1--forms $\omega\in
H^0(K\tilde L)$, where a section of $K\tilde L$ is holomorphic if and
only if it is $d^\nabla$--closed, \cite[Sec.\ 2.3]{Klassiker}.
Generally, the construction of the complex structure of $L$ from the
complex structure of $\tilde L$ requires $\omega$ to have no zeros.
Locally, there always exist holomorphic 1--forms $\omega$ without
zeros.

\begin{theorem} 
\label{t:backward_construction}
Let $\tilde f: M \to \HP^{n-1}$ be a Frenet curve and $\omega\in
H^0(K\tilde L)$ a holomorphic 1--form without zeros.  Then there
exists a Frenet curve $f: M \to \HP^n$ with monodromy whose tangent
curve is $\tilde f$.
\end{theorem}
\begin{proof}
  Since $d^\nabla \omega =0$ there exists a section
  $\tilde\psi\in\Gamma(\tilde V)$ of the trivial $\H^n$--bundle
  $\tilde V$ with translational monodromy such that
 \[
 \nabla \tilde\psi = \omega.
\]
Then $\psi = \tilde\psi \oplus 1 $ is a nowhere vanishing section of
the trivial $\H^{n+1}$--bundle $ V = \tilde V \oplus \trivial{}$.  The
line subbundle $L= \psi\H \subset V$ corresponds to a smooth map $f: M
\to \HP^n$ with loxodromic monodromy.  Since $\omega\in H^0(K\tilde
L)$ is nowhere zero, we define $N: M \to S^2\subset \H$ by
\[ 
*\omega =\tilde S\omega = \omega N
\]
where $\tilde S$ is the canonical  complex structure of
$\tilde V$.  Via the splitting $V = L \oplus \tilde V$ the complex
structure $ J\psi = \psi N $ on $L$ defines the complex structure
\[
\hat S = J \oplus \tilde S
\]
on $V$. By construction $\hat S$ stabilizes the flag $L\subset V_1\subset
\ldots \subset V_{n-1}\subset V$ with $V_k = L \oplus \tilde V_{k-1}$
where $\tilde V_{k}$ is the Frenet flag of $\tilde f$.  Identifying
$\tilde V = V/L$ we get $\tilde \delta_{k-1} =\delta_k$ which implies
\[ 
*\delta_k = \hat S\delta_{k} =  \delta_k \hat S 
\]
for $ k=1,\ldots, n$.  Therefore, to see that $V_k$ is the Frenet flag
 of $L$, it suffices to calculate
\[
*\delta_0 \psi = * \omega =  \tilde S\omega
= \hat S\delta_0 \psi
\]
and
\[
*\delta_0 \psi  =* \omega   =   \omega N 
                = \delta_0 \psi N = \delta_0 J \psi = \delta_0 \hat S\psi \, , 
\]
where we used
$ \delta_0 \psi = \pi_L \nabla\psi = \omega$. 

It remains to show that the canonical complex structure $S$ of $f$
extends smoothly into the Weierstrass points. Since $\hat S$ already
stabilizes the flag $V_k$ the canonical complex structure $S$
decomposes in the splitting $V =L \oplus \tilde V$ into
\begin{equation}
\label{eq:S_splitting}
S = \begin{pmatrix} J & B \\ 0 & \tilde S \end{pmatrix}\,,
\end{equation}
where $B: \tilde V\to L$ is smooth away from the Weierstrass points of
$f$.  Furthermore, the trivial connection $\nabla$ on $V$ decomposes
into
\begin{equation}
\label{eq:nabla_splitting}
\nabla = \begin{pmatrix} \nabla^L & 0 \\ 
                            \delta_0 & \tilde\nabla
            \end{pmatrix},
\end{equation}
where $\delta_0\in\Gamma(K\Hom_+(L, \tilde L))$ is the derivative of
$L$. To check that $S$ is smooth on $M$ it suffices to show that $B$
is smooth on $M$.  But the $(1,0)$--part of $\nabla S$ is given by
\[  
 A  = \frac 14(*\nabla S + S\nabla S) =
       \begin{pmatrix} A_L & \tilde \eta\\ 
                        0 & \tilde A + \frac 12 *\delta_0 B
      \end{pmatrix}\,,
\]
where $-2*A_L =(\nabla^LJ)'$, $-2*\tilde A = (\tilde\nabla \tilde S)' $ and
$\tilde\eta\in\Omega^1(\Hom(\tilde V,L))$.  Since $S$ is the canonical
complex structure of $f$ it follows from (\ref{eq:canonical_S}) that
\begin{equation}
\label{eq:ccs_A}
\tilde A + \frac 12 *\delta_0 B =0
\end{equation}
away from the Weierstrass points. By assumption $\delta_0\psi =
\omega$ is nowhere vanishing so that (\ref{eq:ccs_A}) defines $B$
smoothly also in the Weierstrass points.
\end{proof}

\begin{rem}
  One can always recover the original curve from any of its tangent
  curves: given a Frenet curve $f$ and a hyperplane $H\subset
  \H^{n+1}$ then there is, up to scale, a unique nowhere vanishing
  section $\psi \in \Gamma(L)$ with $\nabla \psi = \delta \psi$,
  namely projections of $\varphi_0\in V= L \oplus H$ onto $L$. Let
  $\tilde f$ be the tangent curve of $f$ with respect to $H$ and let
  $\omega = \delta\psi \in H^0(K\tilde L)$. Choosing the appropriate
  constant of integration we recover $f$.
\end{rem}

\begin{rem}\label{rem:baecklund}
  The enveloping construction is closely related to the B\"acklund
  transformation for Willmore surfaces \cite{coimbra} and its
  generalization \cite{onestep} to holomorphic curves $f: M \to\HP^n$.
  Given such a curve $f$ the B\"acklund transformation \cite{onestep}
  constructs new holomorphic curves $f^\sharp: M \to\HP^k$ by
  integration of $k$ many holomorphic 1--forms in $H^0(KL)$ for $1\le
  k\le \dim H^0(KL)$. In the simplest case of $k=1$, the conformal
  immersion $f^\sharp: M \to\HP^1$ given in affine coordinates by
  $df^\sharp=<\alpha,\omega>$ is a B\"acklund transformation of $f$
  with respect to $\alpha\in (\H^{n+1})^*$.  Here $\omega\in H^0(KL)$
  and $\alpha|_L\in \Gamma(L\invers)$ are assumed to be nowhere
  vanishing sections.
  
  But an envelope $\hat f: M \to\HP^{n+1}$ of $f$ also arises from
  integrating $\omega\in H^0(KL)$. In fact, we see that the nowhere
  vanishing section $\psi\in\Gamma(\hat L)$ constructed in the proof
  of Theorem \ref{t:backward_construction} satisfies
 \[
 d<\hat\alpha,\psi> =  <\hat\alpha,\nabla\psi> = <\alpha,\omega> = df^\sharp\,,
\]
where $\hat\alpha$ is an extension of $\alpha$ to a form in
$(\H^{n+2})^*$.  In particular, this shows that the B\"acklund
transform $f^\sharp$ is given by the projection of the envelope $\hat
f$ onto a suitable $\HP^1\subset\HP^{n+1}$.
\end{rem}

\section{Osculates and envelopes of Willmore spheres}


It is a classical fact that the mean curvature sphere congruence of a
Willmore surface in 3-space is harmonic. The analog of the mean
curvature sphere in the setting of Frenet curves $f: M \to \HP^n$ is
the canonical complex structure.  We will show  that a
Frenet curve is Willmore if and only if the canonical complex
structure is harmonic. Moreover, the constructions of the previous
sections preserve Willmore spheres in $\HP^n$.

As an application, we can construct Willmore spheres in $S^4$ and
minimal spheres with planar ends in $\R^4$ from rational curves in
$\CP^{2n+1}$: the twistor projection \cite[Lemma 2.7]{Klassiker} of a rational
curve in $\CP^{2n+1}$ is a Willmore curve in $\HP^n$. Applying the
tangent construction gives Willmore spheres in $\HP^1 = S^4$ which
themselves are generally not twistor. Therefore, stereographic
projection from an appropriately chosen point on the surface yields a
minimal surface in $\R^4$ with planar ends \cite{montiel}.

\begin{definition}\label{def:Willmore}
  A Frenet curve $f: M \to \HP^n$ is called \emph{Willmore curve} if $f$ is
  critical for the \emph{Willmore energy} 
\begin{equation}
\label{eq:willmore_energy_ccs}
 W = 2 \int_M <A \wedge *A> \,,
\end{equation}
under compactly supported
  variations by Frenet curves, where we also allow the conformal
  structure on $M$ to vary.
  
  Here $A$ is the $(1,0)$--part, (\ref{eq:nabla_S}), of the derivative
  $\nabla S$ of the canonical complex structure $S$.  For an
  endomorphism $B$ we let $<B> := \frac 14 \tr_{\R} B$ be the real
  trace of $B$. Note that $\End_\pm$ are perpendicular with respect to
  this trace inner product.  $S^2=-1$ implies that $A$ anticommutes
  with $S$ so that $A\in\Gamma(K\End_-(V))$.
\end{definition}
In case $f: M \to \R^4$ is an immersion, and thus a Frenet
curve into $\HP^1$, the above definition gives the usual Willmore
energy
\[
W(f) = \frac{1}{2}\int_M |H|^2 - K - K^\perp\,,
\]
where $H$ is the mean curvature of $f$, $K$ the Gaussian curvature and
$K^\perp$ the curvature of the normal bundle of $f$ all computed with
respect to the induced metric on $M$. The critical points of this
functional are called Willmore surfaces and have a long history
attached to them: \cite{blaschke}, \cite{willmore}, \cite{weiner},
\cite{liyau}, \cite{Bryant}, \cite{Ejiri}, \cite{simon},
\cite{montiel}.\kommentar{}
 The next theorem, linking the Willmore condition with
harmonicity, is the natural generalization of the corresponding
theorems in 3 and 4--space.
\begin{theorem}
\label{t:Willmore}
Denote by $\Zz = \{ S\in\End(\H^{n+1}) \mid S^2 =-1\}$ the space of
complex structures on $\H^{n+1}$.
A Frenet curve $f: M \to \HP^n$ is Willmore if and only if the
canonical complex structure $S: M \to \Zz$ is harmonic, that is to
say, $d^\nabla *A =0$.
\end{theorem}
\begin{proof}
  To compute the Euler--Lagrange equation of the Willmore energy, we
  adopt the following perspective: rather then varying the holomorphic
  curve, we vary the background connection by gauge transformations.
  Since $\Gl(n+1,\H)$ acts transitively on the space of flags in
  $\H^{n+1}$ together with complex structures stabilizing the flag,
  these two points of view are equivalent. In particular, $L \subset
  V_t$ is a Frenet curve over $M$ with canonical complex structure $S$
  where we denote by $V_t$ the bundle $V$ equipped with the trivial
  connection $\nabla_t = B_t \nabla B_t \invers$.  Calculating the
  infinitesimal variation of $W$ we get
\begin{equation}
\label{eq:W_dot}
 \dot W =  \int_{M}  4 <\dot A \wedge *A> + 2 < A \wedge \dot * A>\,.
 \end{equation}
 From equation \eqref{eq:nabla_S} we have $ 4 A_t = *_t\nabla_t S +
 S\nabla_t S $ and therefore
\[ 
4 \dot A  = \dot* \nabla S + 2 *\omega_-S + 2 S\omega_-S =
\dot *\nabla S + 4\omega_-'\,,
\]
where $\omega = \dot \nabla = -\nabla \dot B\in\Omega^1(\End(V))$ is
the infinitesimal variation of $\nabla$ and the subscript $\pm$
denotes the projection into $\End_\pm$.  Since $S$ is the canonical
complex structure, we have \eqref{eq:nabla_S} and
\eqref{eq:canonical_S} so that
\begin{eqnarray*}
<\dot * \nabla S \wedge *A> &=& 2 < \dot * *Q\wedge *A> - 2<\dot * *A\wedge *A >\\
&=&- 2<\dot * *A\wedge *A > = 2<*\dot * A\wedge *A > =2<\dot * A \wedge A>\\
& =& -2<A \wedge \dot * A>\,,
\end{eqnarray*}
where we also have used that $\dot *$ and $*$ anti-commute.  Combining
these formulas we get
\[
 \dot W = 4 \int_{M}   <\omega_-'\wedge *A>\,.
\]
Now observe that $*\omega''_- = \omega''_- S$ which implies that
$\omega''_- \wedge *A =0$ by type. Finally, using Stokes Theorem and
recalling that $\End_\pm$ are perpendicular, we obtain
\[
\dot W =  4 \int_{M} <\omega \wedge *A> = 4 \int_{M} <\dot
Bd^\nabla*A>
\]
for any $\dot B\in\Gamma(\End(V))$ with compact support in $M$. In
particular, if $f$ is Willmore, then $d^\nabla *A =0$. By the standard
arguments \cite[Prop.\ 5]{coimbra} this is the Euler--Lagrange equation for $S:
M \to \Zz$ being harmonic. 
\end{proof}

\begin{rem} 
  Any holomorphic curve $f: M \to\HP^n$ has a Willmore energy, namely
  the Willmore energy of the (quaternionic) holomorphic line bundle
  $L\invers$ which coincides with the functional
  (\ref{eq:willmore_energy_ccs}) for Frenet curves \cite[Def.\
  2.5]{Klassiker}.  Thus, a natural definition for a Willmore curve
  would be a holomorphic curve $f: M\to\HP^n$ critical with respect to
  compactly supported variation by holomorphic curves
  $f_t:M_t\to\HP^n$ including variations $M_t$ of the conformal
  structure on $M$. The arguments in Theorem~\ref{t:Willmore} then
  imply that the canonical complex structure $S$ of $f$ is harmonic
  away from the discrete Weierstrass points $D\subset M$ of $f$.
  Unfortunately, this does not guarantee that $S$ extends into the
  Weierstrass points: the example in the beginning of
  Section~\ref{sec:frenet} of twistor projections into $\HP^n$ of
  complex holomorphic curves in $\CP^{2n+1}$ shows that $S$ can become
  singular on $D$.  But if we knew a priori that $S$ is continuous on
  $D$, then the recent regularity result \cite{helein} for harmonic
  maps into pseudo--Riemannian manifolds will show that $S$ is indeed
  smooth in $D$ and hence $f$ by Lemma~\ref{lem:S_smooth} a Frenet
  curve.  In other words, we could replace the assumption of
  Definition~\ref{def:Willmore} that a Willmore curve is Frenet by the
  condition that the canonical complex structure $S$ extends
  continuously into the Weierstrass points and still have
  Theorem~\ref{t:Willmore} valid.
\end{rem}

So far we have seen that the osculating and enveloping constructions
preserve Frenet curves and thus conformality. The situation becomes
much more subtle when considering Willmore curves: it turns out that
the osculates and envelopes of Willmore curves are examples of
constrained Willmore curves. These are Frenet curves critical for the
Willmore energy under compactly supported variations by Frenet curves
preserving the conformal structure of $M$. The theory of constrained
Willmore surfaces in 3--space has only recently been given firm
foundation \cite{christoph_paul_ulrich} and its generalization to
Frenet curves in $\HP^n$ is little understood at present. Therefore,
we restrict ourselves from now on to Willmore spheres in $\HP^n$ in
which case the conformal constraint is void.
\begin{theorem}
\label{t:tangent_willmore}
The tangent curve of a Willmore sphere $f: S^2 \to \HP^n$ is Willmore.
\end{theorem}
\begin{proof} 
Let $\tilde f: S^2 \to P(H)$ be the tangent curve of $f$ with respect
to the hyperplane $H\subset \H^{n+1}$, i.e., $\tilde L = H
\cap V_1$.  Using the splitting $V =L \oplus H$, we decompose the
canonical complex structure $S$ of $f$ into
\begin{equation*}
S = \begin{pmatrix} J & B \\ 0 & \tilde S \end{pmatrix}.
\end{equation*}
Here $J\in\Gamma(\End(L))$ is the complex structure on $L$ and $\tilde
f$ is a Frenet curve with canonical complex structure $\tilde S$.
Moreover, the $(0,1)$--part of $\nabla S$ calculates to
\begin{equation}
\label{eq:Q}
Q = \frac 14(S\nabla S - *\nabla S)
  = \begin{pmatrix} 0 & \eta \\ 0 &\tilde Q\end{pmatrix}
\end{equation}
with $\eta\in\Omega^1(\Hom(H,L))$ and $ 2*\tilde Q = (\tilde\nabla
\tilde S)'' $.  In order to show that $\tilde f$ is Willmore, we have
to calculate $d^{\tilde\nabla}*\tilde A=0$ which by (\ref{eq:nabla_S})
is equivalent to $d^{\tilde\nabla} *\tilde Q =0$.  Since $f$ is
Willmore
\[
 0 =d^\nabla *Q 
  = \begin{pmatrix} 0 &d^{\tilde\nabla, \nabla^L}*\eta\\
                    0&d^{\tilde\nabla}*\tilde Q + \delta \wedge *\eta
     \end{pmatrix}\,.
\]
Flatness of $\nabla$ implies \eqref{eq:nabla_splitting} that $\delta$
is closed and thus there exists $C\in\Gamma(\Hom(L,\tilde L))$ with
\begin{equation*}
\label{eq:C}
\nabla C =\delta.
\end{equation*}

Now, $C*\eta\in\Omega^1(\End(H))$ is a 1--form with values in  
\[
 \tilde\Rr =\{ R\in\End(H) \mid R\tilde V_{n-2} =0, R H \subset \tilde L\}
\]
which satisfies
\begin{equation}
\label{eq:dcomega}
 d^{\tilde \nabla}(C*\eta)= \delta\wedge*\eta + C d^{\nabla^L, \tilde\nabla}*\eta = - d^{\tilde\nabla}*\tilde Q.
\end{equation}
Corollary \ref{cor:constraint_sphere} below now implies that $\tilde f$
is Willmore.
\end{proof}

As an application of this theorem, we construct minimal spheres with
planar ends in $\R^4$:
 \begin{cor}
   Let $h: S^2 \to\CP^{2n+1}$ be a rational curve whose $n^{\text{th}}$
   osculating space $W_n$ does not contain a quaternionic subspace,
   i.e, $W_n \oplus W_n j = \C^{2n+2} = \H^{n+1}$. Then $h$ gives rise
   to a $2(n-1)$--dimensional family of minimal spheres in $\R^4$ with
   planar ends via twistor projection and osculating construction.
\end{cor}
\begin{proof} Under our assumption on the complex holomorphic curve
  $h$ it is shown in \cite[Lemma 2.7]{Klassiker} that the twistor
  projection $f: M \to\HP^n$ of $h$ has the smooth canonical complex
  structure $S$ given by $S|_W = i$ and is Willmore.  For a generic
  choice of complementary hyperplane $H\subset \H^{n+1}$ the tangent
  curve will not be twistor.  Proceeding successively, we obtain a
  $2(n-1)$--dimensional family of Willmore spheres in $\HP^1$ which
  are not twistor projections from $\CP^3$. Therefore, by
  stereographic projections \cite{montiel}, \cite{coimbra}, we obtain
  a $2(n-1)$--dimensional family of minimal spheres with planar ends
  in $\R^4$.
\end{proof}
\begin{figure}[h]
\begin{center}
\epsfig{file=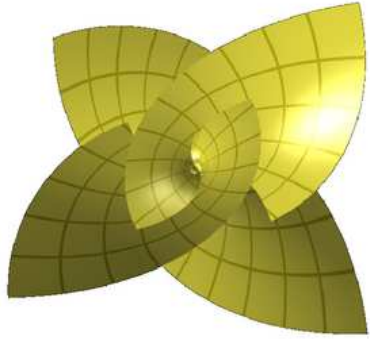,width=0.45\linewidth}
\caption{Minimal sphere in $\R^4$ with planar ends obtained by tangent
  construction on a twistor projection of a holomorphic curve in
  $\CP^5$, \cite{ulisdiss}}
\label{fig:tangent0}
\end{center}
\end{figure}

\begin{rem}
  It is shown in \cite{soliton} that any minimal sphere in $\R^4$ with
  planar ends arises from a rational curve into $\CP^{2n+1}$ via
  twistor projections and osculating constructions.
\end{rem}

As a further application we obtain a splitting of $V = \trivial{n+1}$
into Willmore spheres which is analogous to  the splitting in 
 Corollary~\ref{c:Grothendieck} for Frenet curves:
\begin{cor}
  Let $f: S^2 \to \HP^n$ be a Willmore sphere. Then $V$ splits into a
  direct sum $V = \oplus_{i=0}^n L_i$ of Willmore spheres $f_i :S^2
  \to \HP^{n-i}$ where $f_0 = f$ and $f_n$ is a point in $\HP^1$.
\end{cor}

According to Theorem \ref{t:backward_construction} the construction of
a Frenet curve $f: S^2 \to\HP^n$ with a given first osculating
curve $\tilde f: S^2 \to \HP^{n-1}$ requires the prescription of
tangents $\omega \in H^0(K\tilde L)$. If $\tilde f$ is Willmore there
are natural choices for such $\omega$. Since $d^{\tilde\nabla}*\tilde
A =0$ we can view $\tilde A \in H^0(K\Hom_+(\overline{\tilde V},
\tilde V))$ as a complex holomorphic bundle map. Therefore its kernel
defines a smooth codimension $1$ subbundle $W \subset \tilde V$
provided $\tilde A \not\equiv 0$, i.e., $\tilde f$ is not a twistor
projection. By transversality, we can choose a non--zero $b\in \tilde
V$ such that $\tilde V = W \oplus b\H$. Then
\[
 \omega = *\tilde A b\in H^0(K\tilde L)
\]
is a holomorphic section which vanishes at the zeros of $\tilde A$.
Theorem \ref{t:backward_construction} requires $\omega$ to have no
zeros for the construction of the enveloping Frenet curve $f$.
Nevertheless, the regularity of Willmore surfaces enables us to extend
the construction across the zeros of these specially chosen $\omega$.

\begin{theorem}
\label{t:reverse_Willmore}
Every  Willmore sphere $\tilde f: S^2 \to \HP^{n-1}$, which is not a
twistor projection of a holomorphic curve in $\CP^{2n-1}$,  is a 
tangent curve of a Willmore sphere $f: S^2 \to \HP^n$.
\end{theorem}

\begin{proof} 
  We first prove that $f$, as constructed in Theorem
  \ref{t:backward_construction}, is Frenet. Recall that $f:
  S^2\to\HP^n$ is given as $L =\psi\H$ where $\delta\psi = \omega =
  *\tilde Ab$. The bundle $W \subset \tilde V$ defined by $\ker\tilde A$ is
  stable under the canonical complex structure $\tilde S$ of $\tilde
  f$. Let $\beta\in \Gamma(\tilde V^*)$ with $<\beta, W> =0$ and $
  <\beta, b> =1$.  Then
\[
J \psi = -\psi < \beta,\tilde S b>
\]
defines a complex structure on $L$ and $f$ is admits a Frenet flag
(\ref{eq:frenet_curve0}), (\ref{eq:frenet_curve}) with respect to
the complex structure $J \oplus \tilde S$ by the proof of
Theorem \ref{t:backward_construction}. 

We now show that the canonical complex structure extends smoothly
across the Weierstrass points of $f$.  Away from these points the
canonical complex structure of $f$ can be expressed by
\[ 
S = \begin{pmatrix} J & B \\ 0 & \tilde S \end{pmatrix}
\]
in the splitting $V = L \oplus\tilde V$. Then
\[  
 A  = \frac 14(*\nabla S + S\nabla S) =
       \begin{pmatrix} A_L & \tilde \eta\\ 
                        0 & \tilde A + \frac 12 *\delta B
      \end{pmatrix}
\]
and since $\Im A \subset L$ by \eqref{eq:canonical_S}, we get
\[
 \tilde A + \frac 12 *\delta B =0
\]
away from the Weierstrass points of $f$.
But the smooth bundle map $\psi \beta \in\Gamma(\Hom(\tilde V, L))$
satisfies
\[
\tilde A  = \tilde A b\beta = -*\delta \psi \beta 
\]
and therefore $B = 2\psi\beta$ is smooth across the Weierstrass
points. Thus, the canonical complex structure extends smoothly across
the Weierstrass points of $f$.

It remains to show that $f$ is Willmore. Using \eqref{eq:Q} and the
fact that $\tilde f$ is Willmore, i.e., $d^{\tilde \nabla}*\tilde Q
=0$, we see that
\[
 d^\nabla(*Q -\pr_L *Q)= \begin{pmatrix} 0 & 0 \\ 0 & d^{\tilde\nabla} *\tilde Q \end{pmatrix} = 0.
\]
As in the proof of the previous theorem,  we denote by
\[
\Rr = \{ R\in \End(V) \mid RV_{n-1} =0, RV \subset L\}\,,
\]
where $V_{n-1}$ is the $n-1^{\text{st}}$ osculating bundle of $f$.
Then $\pr_L*Q$ takes values in $\Rr$ and Corollary
\ref{cor:constraint_sphere} below implies that $f$ is Willmore.
\end{proof}

\begin{rem}
 From Remark \ref{rem:baecklund} and Theorem \ref{t:reverse_Willmore}
we see that the enveloping construction for Willmore spheres coincides
up to a suitable projection with the B\"acklund transformation for
Willmore spheres.
\end{rem}

\begin{figure}[h]
\begin{center}
\epsfig{file=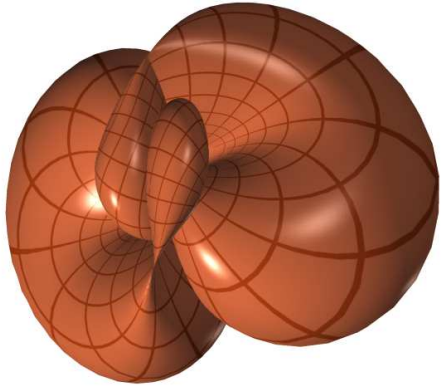,width=0.45\linewidth}
\caption{Constrained Willmore torus in $\R^4$ obtained by tangent
  construction on a twistor projection of an elliptic curve in $\CP^5$,
  \cite{ulisdiss}}
\label{fig:tangent2}
\end{center}
\end{figure}

\section{Technical Lemmas}
We conclude our paper by providing various technical lemmas and the
corollary used in the proofs of Theorem \ref{t:tangent_willmore} and
Theorem \ref{t:reverse_Willmore}.  Given a Frenet curve $f: M \to
\HP^n$ with osculating flag $V_k$, we let
\[ 
\Rr = \{ R\in \End(V) \mid RV_{n-1} =0, \ RV \subset L\}=\Hom(V/V_{n-1},L)\,.
\]

\begin{lemma}
\label{l:Q}
If $f: M \to \HP^n$ is a Frenet curve with canonical complex structure
$S$ then
\[
 d^\nabla *Q \in \Omega^2(\Rr)
\]
where $2*Q = (\nabla S)''$.
\end{lemma}
\begin{proof}

 Recall that by \eqref{eq:nabla_S} we have $d^\nabla*A = d^\nabla *Q$.
 Since
  $A\in\Gamma(K\Hom_-(V,L))$ we obtain
\[ 
\pi_L d^\nabla *Q = \pi_L d^\nabla*A = \delta_0\wedge *A =0\,,
\]
where $\pi_L: V \to V/L$ and $\delta_0 =
\pi_L\nabla|_L\in\Gamma(K\Hom_+(L,V/L))$ is the derivative of $f$.
This shows that $d^\nabla* Q$ is $L$--valued.

On the other hand, for a section $\psi\in\Gamma(V_{n-1})$, we compute
\[
 (d^\nabla*Q)\psi = -*Q\wedge \nabla\psi = -*Q\wedge \delta_{n-1}\psi =0\,.
\]
Here we used that $QV_{n-1} =0$, $\delta_{n-1}=
\pi_{V_{n-1}}\nabla|_{V_{n-1}}\in\Gamma(K\Hom_+(V_{n-1},V/V_{n-1}))$
and  $*Q = QS$. In other words, $d^\nabla*Q$ vanishes on $V_{n-1}$,
and therefore $d^\nabla*Q\in\Omega^2(\Rr)$.
\end{proof}

\begin{lemma}
\label{l:tuy}
 Let $f:M \to \HP^n$ be a Frenet curve  with  complex structure $S$.
Then  $\eta\in \Omega^1(\Rr)$ and $d^\nabla\eta\in\Omega^2(\Rr)$ imply 
$
 \eta\in\Gamma(K\Rr_+).
$
\end{lemma}
\begin{proof}
Since $\eta$ is $L$--valued,
 \[
0 = \pi_Ld^\nabla\eta =\delta_0\wedge \eta
\]
implies $\eta\in\Gamma(K\Rr)$.  Furthermore, for
$\psi\in\Gamma(V_{n-1})$
\[
0=(d^\nabla\eta)\psi = -\eta\wedge \delta_{n-1}\psi
\] 
shows that $*\eta=\eta S$ and thus $\eta\in\Gamma(K\Rr_+)$.
\end{proof}

Let $f: M \to \HP^n$ be a Frenet curve with complex structure $S$ on
the trivial $\H^{n+1}$--bundle $V$. From \eqref{eq:nabla_S} it follows
that the $S$--anticommuting part of $\nabla$ is given by $ A + Q$.
The $S$--commuting part of $\nabla$ is a complex connection $ \nabla_+
= \partial + \delbar$ whose $(0,1)$--part $\delbar$ defines a complex
holomorphic structure on $V$. From \eqref{eq:frenet_curve} it follows
that 
\[
0 = (\delta_k)''_+ = \pi_{V_k} \delbar|_{V_k}
\]
so that $V_k\subset V$ is $\delbar$--stable.  Therefore $\delbar$
induces a complex holomorphic structure on the complex line bundle
$K\Rr_+$.

\begin{lemma}
\label{l:etaholomorph}
Let $f: M \to \HP^n$ be a Frenet curve with canonical complex
structure $S$ and $\eta\in\Omega^1(\Rr)$ with $ d^\nabla \eta =
d^\nabla *Q$. Then $ \eta \in H^0(K\Rr_+)$ is a complex holomorphic
section.
\end{lemma}
\begin{proof}
From Lemma \ref{l:Q} and Lemma \ref{l:tuy} it follows 
 that $\eta\in\Gamma(K\Rr_+)$.  Since 
 $\nabla = \nabla_+ + A + Q$ we obtain
\[
 d^\nabla\eta = d^{\nabla_+} \eta + [A\wedge\eta] + [Q\wedge\eta]\,.
\]
But $\eta\in\Gamma(K\Rr_+)$, and $A$ and $ Q$ anticommute with $S$ so that
\[
 (d^\nabla\eta)_+ = \dbar \eta\,.
\]
On the other hand, 
\[
 d^\nabla *Q = d^{\nabla_+}*Q + [A\wedge*Q] + [Q\wedge *Q] = d^{\nabla_+}*Q\,,
\]
where we used that $[A\wedge*Q] = 0$ and $[Q\wedge *Q]=0$ by type and
symmetry considerations. This shows $(d^\nabla*Q)_+=0$ and therefore
\[
\delbar\eta= (d^\nabla\eta)_+ = (d^\nabla*Q)_+=0\,,
\]
which proves that $\eta\in H^0(K\Rr_+)$.
\end{proof}

\begin{lemma}
\label{l:degree} If $f: M \to \HP^n$ is a Frenet curve with complex
structure $S$ then
\[ 
\delta_{n-1} \circ \ldots \circ \delta_0 \in H^0(K^n\Hom_+(L,V/V_{n-1}))
\]
is complex holomorphic and hence the degree of the complex line bundle
$K\Rr_+$ is
\[ 
\deg K\Rr_+ = (n+1)\deg K - \ord(\delta_{n-1} \circ \ldots \circ \delta_0).
\]
In particular, if $M=S^2$  then $\deg K\Rr_+ <0$, so that $K\Rr_+$ has
no global holomorphic sections.
\end{lemma}
\begin{proof}
  From \eqref{eq:frenet_curve} and  \cite{Le}   it follows that the
  derivative $\delta_k$ of the $k^{\text{th}}$ osculating flag $V_k$
  is a complex holomorphic section of the complex line bundle
  $K\Hom_+(V_k/V_{k-1},V_{k+1}/V_k)$.  Therefore
\[ 
\delta_{n-1} \circ \ldots \circ \delta_0 \in H^0(K^n\Hom_+(L,V/V_{n-1}))
\]
and the degree of the complex line bundle $K^n\Hom_+(L,V/V_{n-1})$
is given by
\[ 
\ord(\delta_{n-1} \circ \ldots \circ \delta_0) 
  =
 \deg(K^n\Hom_+(L,V/V_{n-1}))\,.
\]
Since the degree of a complex quaternionic line bundle is  the
degree of the underlying complex line bundle, \cite{Klassiker},
we have
\[
 \deg(K^n\Hom_+(L,V/V_{n-1}))
= n\deg K +\deg (V/V_{n-1}) -\deg L\,.
\]
Recalling that $\Rr_+=\Hom_+(V/V_{n-1},L)$, we finally obtain
\[
 \deg K\Rr_+ 
=
 \deg K + \deg L - \deg(V/V_{n-1})
 = 
(n+1)\deg K - \ord(\delta_{n-1} \circ \ldots \circ \delta_0).
\]
\end{proof}

Combining the last two lemmas and Theorem \ref{t:Willmore} yields the
following
\begin{cor}
\label{cor:constraint_sphere}
Let $f: S^2 \to \HP^n$ be a Frenet curve with canonical complex
structure $S$ satisfying
\[
d^\nabla(*Q + \eta) =0 
\]
for some $\eta\in\Omega^1(\Rr)$. Then $f$ is Willmore, i.e.,
$d^\nabla*Q =0$.
\end{cor}



\bibliographystyle{alpha}

\bibliography{doc,doc_dpw}

\end{document}